\documentclass[10pt]{article}
\usepackage[tbtags]{amsmath}
\usepackage{epsfig,amstext,amssymb,amsthm,latexsym}
\allowdisplaybreaks[4]

\catcode`\@=11 \@addtoreset{equation}{section} \catcode`\@=12

\newcommand{\nwc}{\newcommand}
\nwc{\COM}[1]{}

\nwc{\vs}[1]{\vskip #1 cm}

\newtheorem{theo}{Theorem}[section]
\newtheorem{sat}[theo]{Proposition}
\newtheorem{de}[theo]{Definition}
\newtheorem{lem}[theo]{Lemma}

\newtheorem{korr}[theo]{Corollary}
\newtheorem{remark}[theo]{Remark}

\newcommand{\netheo}[1]{{Theorem \ref{#1}}}

\newcommand{\kb}[1]{\boldsymbol{#1}}
\newcommand{\vk}[1]{\kb{#1}}

\def\FRE{\mbox{Fr\'{e}chet }}

\def\x{\vk{x}}

\def\X{\vk{X}}

\newcommand{\ve}{\varepsilon}
\newcommand{\abs}[1]{\lvert #1 \rvert}
\newcommand{\Abs}[1]{ \Bigl \lvert #1 \Bigr \rvert}

\newcommand{\norm}[1]{\lVert #1 \rVert}

\newcommand{\E}[1]{\mbox{\rm$\vk{E}$}\{#1\}}

\newcommand{\pk}[1]{\mbox{\rm$\vk{P}$} \{#1\} }

\newcommand{\pb}[1]{\mbox{\rm$\vk{P}$}\Bigl \{#1 \Bigr \}}

\newcommand{\R}{\!I\!\!R}

\newcommand{\inr}{\in \R}

\newcommand{\ldot}{,\ldots,}

\newcommand{\limit}[1]{\lim_{#1 \to   \infty}}

\newcommand{\as}{ \stackrel{a.s}{\to}}
\newcommand{\todis}{\stackrel{d}{\to}}
\newcommand{\toprob}{ \stackrel{p}{\to}}

\newcommand{\equaldis}{\stackrel{d}{=}}

\newcommand{\BQN}{\begin{eqnarray}}
\newcommand{\EQN}{\end{eqnarray}}
\newcommand{\BQNY}{\begin{eqnarray*}}
\newcommand{\EQNY}{\end{eqnarray*}}
\newcommand{\BS}{\begin{sat}}
\newcommand{\ES}{\end{sat}}
\newcommand{\BL}{\begin{lem}}
\newcommand{\EL}{\end{lem}}
\newcommand{\BT}{\begin{theo}}
\newcommand{\ET}{\end{theo}}
\newcommand{\BK}{\begin{korr}}
\newcommand{\EK}{\end{korr}}
\newcommand{\BD}{\begin{de}}
\newcommand{\ED}{\end{de}}
\newcommand{\BIT}{\begin{itemize}}
\newcommand{\EIT}{\end{itemize}}
\newcommand{\BDI}{\begin{description}}
\newcommand{\EDI}{\end{description}}

\newcommand{\QED}{\hfill $\Box$}
\newcommand{\IF}{\infty}

\def\fracl#1#2{\biggr( \frac{#1}{#2} \biggl) }

\newcommand{\prooftheo}[1]{ \textsc{Proof of Theorem} \ref{#1} }

\newcommand{\proofkorr}[1]{\textsc{Proof of Corollary} \ref{#1}}

\def\X{\vk{X}}
\def\Y{\vk{Y}}
\def\SI{\Sigma}

\def\LM{\lambda_{(m)}}
\def\U{\vk{U}}
\def\Um{\vk{U}^{(m)}}
\def\Umd{\vk{U}^{(d-m)}}
\def\HL{H\"usler et al.\ (2002)}

\begin{document}

\centerline{\Large Asymptotics of the Norm of  }
        \vskip .4 cm
\centerline{\Large Elliptical Random Vectors}

        \vskip 0.8 cm
        \centerline{\large Enkelejd Hashorva}
        \vskip 0.8 cm

        \centerline{\textsl{Department of Statistics, University of Bern}}
        \centerline{\textsl{Sidlerstrasse 5, CH-3012, Bern, Switzerland}}
        \centerline{\textsl{Email: enkelejd.hashorva@stat.unibe.ch}}

\today{}
        \vskip 1.4 cm

{\bf Abstract:} In this paper we consider elliptical random vectors
$\X$ in $\R^d,d\ge 2$ with stochastic representation $  A R \U,$ where
$R$ is a positive random radius independent of the random vector
$\U$ which is uniformly distributed on the unit sphere of $\R^d$
and $A\in \R^{d\times d}$ is a given matrix. Denote by $\norm{\cdot }$ the Euclidean norm in $\R^d$, and let $F$ be the distribution function of $R$.  The main result of this paper is an asymptotic expansion of the probability
$\pk{\norm{\X} > u}$ for $F$ in the Gumbel or the Weibull
max-domain of attraction. In the special case that $\X$ is a mean zero Gaussian random vector our result coincides with the one derived in \HL.

{\it Key words and phrases:  Elliptical distribution;  Gaussian distribution; Kotz Type distribution;
Gumbel max-domain of Attraction; Tail approximation; Density convergence; Weak convergence.}

\section{Introduction}
Let $\X$ be a  mean zero Gaussian random vector in $\R^d, d\ge 2$
with underlying covariance matrix $\SI$.   If $\LM:=\lambda_1\ge \lambda_2 \ge \cdots \ge \lambda_d$
are the ordered eigenvalues of the matrix $\SI$, then in the light of Theorem 1 in \HL \, we have the
asymptotic expansion  (set $C_*:=\prod_{j=m+1}^d(1 - \lambda_j/\LM)^{-1/2}  $ and $C_*:=1$ if $m=d$)
\BQN\label{eq:huesler}
\pb{\norm{\X} > \sqrt{u\LM} }&=& C_*\frac{2^{1- m/2}}{\Gamma(m/2)} u^{m/2- 1}\exp(-u/2), \quad u\to \IF,
\EQN
where $\norm{\vk{x}}$ stands for the Euclidean norm of $\vk{x}\inr^d$, $m$ is the multiplicity of $\lambda_1$ i.e., 
$m:=\# \{j: \lambda_j=\lambda_1, 1\le j\le d \}$, and
 $\Gamma(\cdot)$ is the Gamma function.\\
It is well-known (see e.g., Cambanis et al.\ (1981), or Fang et al.\ (1990)) that
$\X$ possesses the stochastic representation
\BQN\label{eq:stoch:e}
\X&\equaldis&  A R\U,
\EQN
with $R>0$ such that $R^2$ is Chi-squared distributed with $d$ degrees of freedom
being further independent of $\U$ which is uniformly distributed on the
unit sphere of $\R^d$ and $A$ is a $d\times d$ real matrix satisfying $A A^\top= \SI$. Here $\equaldis$ and $\top$
stand for equality of distribution functions and the transpose sign, respectively.

If we drop the distributional assumption on $R$ assuming simply that $R>0$ almost surely with some unknown distribution function $F$ with upper endpoint $x_F \in (0,\IF]$, then the random vector $\X$ with stochastic representation \eqref{eq:stoch:e} is an elliptical random vector (see Cambanis et al.\ (1981)). \\
In this paper we focus our interest in possible generalisation of \eqref{eq:huesler} considering some
general elliptical random vector $\X$. It is clear that the  asymptotics in \eqref{eq:huesler} is related to the tail asymptotics of $F$. In the special case that $R^2$ is Chi-squared distributed with $d$ degrees of freedom we have
\BQN \label{eq:1:chi}
1- F(u)&=& (1+o(1)) \frac{u^{d-2} \exp(- u^2/2)}{ 2^{d/2 -1}\Gamma(d/2)}, \quad u\to \IF
\EQN
implying that  $F$ is in the max-domain of attraction of the unit Gumbel distribution function $\Lambda$ i.e.,
\BQN \label{eq:rdfd} \lim_{u \uparrow x_F} \frac{1 -F(u+x/w(u))}{1- F(u)} &=& \exp(-x),\quad \forall x\inr,
\EQN
where $w(u)= u, u>0$.  \\
Under the max-domain of attraction assumption \eqref{eq:rdfd} on $F$ we generalise \eqref{eq:huesler} in our main result below (see Theorem 1) to the following asymptotic expansion
\BQN\label{eq:theo:A:1}
\pk{\norm{\X} > u \sqrt{ \LM} } &=&  (1+o(1)) C_*\frac{\Gamma(d/2)}{\Gamma(m/2)}  \fracl{2}{u w(u )}^{(d-m)/2}[1- F(u)],
\quad u \uparrow  x_F
\EQN
and obtain also an asymptotic expansion of the density function of $\norm{\X}$. We derive similar asymptotic results when the distribution function $F$ is in the Weibull max-domain of attraction.
Further, considering an elliptical random sample we apply our asymptotic expansions in order to derive convergence in distribution, density convergence and almost sure convergence of the maximum norms.

Brief outline of the rest of the paper: We proceed with a short section dedicated to our notation and some preliminary results. The main results are presented in Section 3. The proofs of all the results are relegated to Section 4.

\section{Preliminaries}
In this section we introduce several notation and briefly discuss elliptical random vectors and max-domain of attractions.
Given a random variable $R$ with distribution function $F$ (write $R \sim F$) we denote by  $\overline{F}$ the survivor function.  If $X$ is Beta distributed with positive parameters $a,b$, then we denote this by $X \sim Beta(a,b)$.
The Beta distribution $Beta(a,b)$ possesses the density function $x^{a-1}(1- x)^{b-1}\Gamma(a+b)/(\Gamma(a)\Gamma(b)), x\in (0,1).$\\

If $I$ is a non-empty subset of $\{1 \ldot d\}, d\ge 2,$ then  $\abs{I}$ denotes the number of its  elements.
For any vector $\x=(x_1 \ldot x_d)^\top \inr^d$ we define the subvector $\x_I$ with respect to  $I$
by $\x_I:=(x_i,i\in I)^\top$. The Euclidean norm of $\x_I\inr^m$ is $\norm{\x_I}:= \sqrt{\sum_{i\in I} x_i^2}$.

In the sequel we write  $\U^{(k)}:=(U_1 \ldot U_k)^\top , k=1 \ldot d$ for a random vector uniformly distributed on
the unit sphere of $\R^k$. For notational simplicity we write $\U$ instead of $\U^{(d)}$.

Basic distributional results for spherical random vectors are derived in Cambanis et al.\ (1981).
By Lemma 2 therein for two non-empty disjoint sets $I,J$ such that $I \cup J=\{1 \ldot d\}$,
and a uniformly distributed random vector  $\U$ we have the stochastic representation
\BQN\label{rep:elli0}
\U_I\equaldis \sqrt{W_{m,d}} \Um, \quad \U_J\equaldis (1- W_{m,d})^{1/2} \Umd,
 \EQN
with $W_{m,d} \sim Beta(m/2, (d-m)/2)$.  Furthermore, $\Um,\Umd,W_{m,d}$ are mutually independent. For our investigation it is crucial that the random vector $\U$ is distributional invariant with respect to orthogonal transformations, i.e.,
\BQN\label{eq:otho}
D \vk{U} &\equaldis &\vk{U}
\EQN
for any orthogonal matrix $D\in \R^{d\times d}$. Furthermore, for two square matrices $A,B$ 
\BQN\label{eq:si}
A \vk{U}&\equaldis & B \vk{U}
\EQN
whenever $BB^\top = AA^\top  \inr^{d\times d}$ is valid.

Next we mention some facts from the univariate extreme value theory:
A distribution function $F$ is said to belong to the max-domain of attraction of a univariate
extreme value distribution function $H$, if for constants $a_n>0,b_n,n\ge 1$
\BQN \label{eq:LL}
 \limit{n}\sup_{x\inr} \Abs{F^n(a_nx+ b_n)- H(x)}&=& 0.
 \EQN
Only three choices for $H$ are possible (see e.g., Galambos (1987), Reiss (1989), Embrechts et al.\ (1997), Falk et al.\ (2004),
Kotz and Nadarajah (2005), De Haan and Ferreira (2006), or Resnick (2008)), namely
the \FRE distribution, the Gumbel distribution, or the Weibull distribution.\\
It is well-known that  if  $F$ with upper endpoint $x_F$ is in the max-domain of attraction of the \FRE distribution   $\Phi_\gamma(x)=\exp(-
x^{-\gamma}),x>0, \gamma\in (0,\IF)$, then necessarily  $x_F=\IF$, and furthermore
\BQN\label{eq:PhiF}
 \limit{u} \frac{\overline{F}(xu)}{\overline{F}(u)}&=& x^{-\gamma}, \quad \forall x>0.
 \EQN
If $H= \Lambda$ with $\Lambda(x)=\exp(-\exp(-x)),x\inr$, then as mentioned in the Introduction \eqref{eq:LL} is equivalent to
\eqref{eq:rdfd}. Alternatively we write $F \in MDA(\Lambda,w)$ whenever \eqref{eq:rdfd} is satisfied. \\

When \eqref{eq:LL} holds with $H$ the Weibull distribution function  $\Psi_\gamma(x)= \exp(-\abs{x}^\gamma), x< 0,\gamma\in (0,\IF)$, then
necessarily $x_F$ is finite, and furthermore
\BQN\label{eq:Psi}
 \limit{u} \frac{\overline{F}(x_F- x/u)}{\overline{F}(x_F-1/u)}& =& x^{\gamma}, \quad \forall x>0.
\EQN
Clearly, \eqref{eq:LL} means that the sample maxima converges in distribution (after normalisation). If $F$ possesses a density function $f$, then the convergence in \eqref{eq:LL} can be strengthened in several instances  to local uniform convergence of the corresponding density functions. In the Gumbel case local uniform convergence follows if
for some positive scaling function $w$
\BQN\label{eq:de:G}
\lim_{ u \uparrow x_F}\frac{f(u+ x/w(u))}{f(u)}&=& \exp(-x), \quad \forall x\inr,
\EQN
whereas for $F$ in the Weibull max-domain of attraction it suffices that
\BQN\label{eq:de:W}
\lim_{u\downarrow 0} \frac{ u f(x_F- u)}{\overline{F}(x_F- u)}&=& \gamma.
\EQN
\COM{
and
\BQN
\limit{u} \frac{ u f(u)}{\overline{F}(u)}&=& \gamma,
\EQN
}
See e.g., Reiss and Drees (1992) or Resnick (2008) for more details.

\section{Main Results}
Let $\X=(X_1 \ldot X_d)^\top $ be an elliptical random vector in $\R^d$ with
stochastic representation \eqref{eq:stoch:e},  where $R\sim F, F(0)=0$,  and $A\in \R^{d\times d}$ is a given square matrix.
We set throughout this paper $\SI := A A^\top$ and denote its eigenvalues by
$$\lambda_1 \ge \lambda_2 \ge \cdots \ge \lambda_d\ge 0$$
which exist since the matrix $\SI$ is semi-positive definite. In the following $m$ stands for the multiplicity of the largest eigenvalue $\lambda_1=:\LM$.\\
As already shown in \HL \,  for the Gaussian setup the tail asymptotics of interest
  depends only on  the eigenvalues of $\SI$, but not on the covariance matrix $\SI$ itself. This must be the case also when $\X$ is an elliptical random vector with stochastic representation \eqref{eq:stoch:e}.
Indeed, first note that by \eqref{eq:si} the matrix $\SI$ specifies the distribution of $\X$ (and not the matrix $A$).
Furthermore  \eqref{eq:otho} implies
$$ \norm{\X}^2\equaldis  R^2\biggl( \lambda_1 U_1^2 +   \cdots +\lambda_d U_d^2\biggr),$$
hence by \eqref{rep:elli0}
$$ \norm{\X}^2\equaldis  \LM R^2\biggl( W_{m,d}+  [V_{1}^2\frac{\lambda_{m+1}}{\LM} + \cdots + V_{d-m}^2\frac{\lambda_{d}}{\LM}] (1- W_{m,d})\biggr)$$
is valid with $\vk{V}_{d-m}\equaldis \Umd$ and $W_{m,d} \sim Beta(m/2, (d-m)/2)$. Furthermore, $R, W_{m,d}, \vk{V}_{d-m}$ are mutually independent. Clearly, for any $u>0$
\BQN\label{eq:bounds}
 \overline{F}(\sqrt{u/\LM})= \pk{\LM R^2 > u}  \ge \pk{\norm{\X}^2 \ge u} \ge \pk{\LM R^2 W_{m,d}> u}.
 \EQN
We show next that for $F$ in the Gumbel or the Weibull max-domain of attraction the asymptotic behaviour of $\pk{\norm{\X}^2 > u}$ is similar to that of  $\pk{\LM R^2 W_{m,d}> u}$. If $F$ is in the max-domain of attraction of $\Phi_\gamma, \gamma \in (0,\IF)$, which is equivalent with $X_1$ has distribution function in the max-domain of attraction of $\Phi_\gamma$ (see Hashorva (2006,2007a)), it follows easily that the random vector $(X_1^2 \ldot X_d^2)^\top$ is regularly varying with index $\gamma/2$. Hence the asymptotic behaviour
  of $\pk{\norm{\X}> u}, u\uparrow \IF $ can be easily determined. In particular $\norm{\X}$ has distribution function in the max-domain of attraction of $\Phi_{\gamma}$.   In the sequel we discuss therefore only the Gumbel and the Weibull cases.

\subsection{Asymptotics in the Gumbel Model}
In this section we assume that the distribution function $F$ of the associated random radius $R$
is in the Gumbel max-domain of attraction. Simple instances of distribution functions $F$ in the Gumbel max-domain of attractions are univariate distributions with exponential tails. In our investigation the scaling function $w$ (see \eqref{eq:rdfd}) plays a crucial role. The following asymptotic properties of $w$ are well-known (see e.g., Resnick (2008))
\BQN\label{eq:uu}
 \lim_{u\uparrow x_F} u w(u)=\infty, \quad \text{and  }
\lim_{u\uparrow x_F} w(u)(x_F - u) =\infty \quad \text{if  } x_F< \infty.
\EQN
Furthermore, $w$ can be defined asymptotically via the mean excess function (see e.g., Embrechts et al.\ (1997)) as
\BQN \label{eq:reco1}
w(u)&=& \frac{1+o(1)}{ \E{R-u \lvert R>u}}, \quad u\uparrow x_F.
\EQN
Throughout the rest of the paper  $x_F\in (0,\IF]$ denotes  the upper endpoint of $F$, 
and $C_*$ is a positive constant (also defined in the Introduction)  given by
\BQN \label{cy}
C_*&:=&\prod_{j=m+1}^d(1 - \lambda_j/\LM)^{-1/2},
\EQN
where $C_*:=1$ if $m=d$.\\
Since $\X/x_F,$ with $x_F\in (0,\IF)$ is again an elliptical random vector, and moreover $F(s/x_F),s\inr$ is in the Gumbel or Weibull max-domain of attraction if $F$ is in the Gumbel or Weibull max-domain of attraction, respectively,
 we assume in the following without loss of generality that $x_F= 1$ or $x_F=\IF$. Next, we state the main result for the Gumbel setup.\\

\BT \label{theo:A} Let $\X=A R \U$ be an elliptical random vector in $\R^d,d\ge 2$,
with $A$ a $d\times d$ real matrix, $R \sim F$ a positive random variable independent of $\U$.
 If $F(0)=0, x_F\in \{1, \IF\}$ and $F\in MDA(\Lambda, w)$ with $w$ some positive scaling function,
then we have
\BQN\label{eq:theo:A:1}
\pk{\norm{\X} > u \sqrt{ \LM} } 
&=& (1+o(1)) C_*\frac{\Gamma(d/2)}{\Gamma(m/2)}
 \fracl{2}{u w(u )}^{(d-m)/2}\overline{F}(u), \quad u\uparrow x_F.
\EQN
Furthermore, $\norm{\X}$ possesses the positive density function $h$ with asymptotic behaviour
\BQN\label{eq:theo:A:2}
h(u)  &=& (1+o(1)) w(u) \pk{\norm{\X} > u}, \quad u\uparrow x_F.
\EQN
\ET

We have now the following result:
\BK \label{korr:1} Let $\X, \X^{(1)} \ldot \X^{(n)},n\ge 1$ be independent elliptical random vectors in $\R^d,d\ge 2$ with common distribution function $G$ such that $\X$ satisfies the assumptions of \netheo{theo:A}. Assume for simplicity that $\LM=1$.
Then we have the convergence in distribution
\BQN\label{eq:conv:dis}
\frac{\max_{1 \le j \le n} \norm{\X^{(j)}} - a_n}{b_n}  &\todis & Y \sim \Lambda, \quad n\to \IF,
\EQN
where $$ a_n:= 1/w(b_n), \quad b_n:= H^{-1}(1- 1/n), \quad n> 1,$$
with  $H^{-1}$ the generalised inverse of the distribution function of $\norm{\X}$.\\
Furthermore, \eqref{eq:conv:dis} can be strengthened to the local uniform convergence of the corresponding density functions.
\EK

\bigskip

\begin{remark}
1. \  The scaling function $w$ is self-neglecting (see e.g., Reiss (1989) or Resnick (2008)) i.e.,
\BQN\label{eq:self}
 \limit{n} w(u + x/w(u))/w(u)&=& 1, \quad u\uparrow x_F
\EQN
uniformly for $x$ in compact sets of $\R$. In view of \eqref{eq:uu} under the assumptions of \netheo{theo:A}
\BQNY
\lim_{u\uparrow x_F}\frac{\pk{\norm{\X} > u \sqrt{ \LM} }}{\overline{F}(u)}&=&0.
\EQNY
Hence the upper bound for $\pk{\norm{\X} > u}$ in \eqref{eq:bounds} is not
accurate for $u$ close to $x_F$. However, the lower bound therein turns out to be asymptotically accurate. \\
2. \  Clearly, the local uniform convergence of the density functions of the sample maxima implies the convergence in distribution stated in \eqref{eq:conv:dis}. See Resnick (2008) for deeper results on the density convergence of the univariate sample extremes.\\
\end{remark}

We provide next three illustrating examples:

\bigskip
{\bf Example 1.} (Gaussian random vectors) Let $\X \inr^d,d\ge 2$ be Gaussian random vector with covariance matrix $\Sigma= A A^\top$ and mean zero. As mentioned in the Introduction $\X$ possesses  the stochastic representation \eqref{eq:stoch:e}, and furthermore $R^2$ is Chi-squared distributed with $d$ degrees of freedom. Since $F\in MDA(\Lambda, w)$, with the scaling function $w(u)= u, u>0$ (recall \eqref{eq:1:chi})
\netheo{eq:theo:A:1} yields \eqref{eq:huesler} 
which is obtained in Theorem 1 of \HL.

{\bf Example 2.} ($F$ with finite upper endpoint) Let $\X \equaldis R A \U, R \sim F $ be an elliptical random vector in $\R^d,d\ge 2$. Assume that $F$ has upper endpoint 1,  and furthermore
$$ \overline{F}(u)= (1+o(1))c_1 \exp(- c_2/(1- u)) , \quad u \uparrow 1,$$
with $c_1,c_2$ two positive constants. Since for $w(u)= c_2/(1- u)^2, u\in (0,1)$ and any $s\inr$ we have 
$$ \frac{ \overline{F}(u+ s/w(u))}{ \overline{F}(u)}= (1+o(1))\exp( - c_2 [ 1/(1- u+ s/w(u)) - 1/(1- u)]) \to \exp(-s), \quad u \uparrow 1,$$
then $F\in MDA(\Lambda, w)$. For this example 
\BQNY
\pk{\norm{\X} > u\sqrt{ \LM} } &=&
 (1+o(1)) c_1 C_* \frac{\Gamma(d/2)}{\Gamma(m/2)}
(2 c_2 /(1-u) )^{(d- m)/2} \exp(- c_2/(1-u)), \quad u \uparrow 1.
\EQNY

\bigskip
{\bf Example 3.} (Kotz Type III elliptical random vectors) Let $\X=  A R \U $ be a $d$ dimensional elliptical random vector. Assume that the survivor function $\overline{F}$  of $R$ satisfies
\BQN\label{eq:kotz:Fk}
\overline{F}(u)&=&  (1+o(1))c  u^{N}\exp(-\delta u^\tau ), \quad  u\to \IF,
 \EQN
where $c,\delta,\tau$ are given positive constants and $N\inr$. We refer to such $\X$ as a Kotz Type III elliptical random vector. It follows easily that $F\in MDA(\Lambda,w)$ with the scaling function $w$ given by
\BQN\label{eq:sc:III}
w(u)&:=& \delta \gamma  u ^{\tau -1},\quad u>0.
\EQN
In view of \netheo{eq:theo:A:1} we may thus write
\BQN\label{tail:K:III}
\pk{\norm{\X} > u \sqrt{\LM} }
%
&=&  (1+o(1))
C_*\frac{\Gamma(d/2)}{\Gamma(m/2)}
 ( 2/( \delta \tau ) )^{(d-m)/2}   u ^{\tau(m-d)/2 +N}\exp(-\delta u^\tau), \quad u\to \IF\notag\\
 &=:&  (1+o(1))Ku ^{\alpha}\exp(-\delta u^\tau ), \quad u\to \IF.
 \EQN
Let  $\X^{(1)} \ldot \X^{(n)}, n\ge 1$ be independent random vectors in $\R^d$ with the same distribution function as $\X$, and assume for simplicity that $\LM=1$. Then the convergence in distribution 
\BQNY
 (\max_{1 \le j  \le n} \norm{\X^{(j)}} - a_n)/ b_n &\todis &  Y \sim \Lambda  ,\quad n\to \IF
 \EQNY
holds with constants $a_n,b_n$ defined by
$$ a_n:= (\delta^{-1} \ln n )^{1/\tau -1}/(\delta\tau),
\quad b_n:= (\delta^{-1}\ln n )^{1/\tau}+a_n\Bigl[ \alpha \ln (\delta^{-1} \ln n )/\tau+ \ln K\Bigr], \quad n> 1.
$$
The above convergence in distribution implies the convergence in probability
\BQNY
\frac{\max_{1 \le j  \le n} \norm{\X^{(j)}}}{(\ln n)^{1/\tau}} &\toprob & \delta^{-1/\tau},\quad n\to \IF.
\EQNY
By the Barndorff-Nielsen  criterion for the almost sure stability of the sample maxima (see Barndorff-Nielsen (1963),
Resnick and Tomkins (1973) or Tomkins (1986))  we retrieve further the almost sure convergence
\BQN\label{eq:asT}
\frac{\max_{1 \le j  \le n} \norm{\X^{(j)}}}{(\ln n)^{1/\tau}} &\as  & \delta^{-1/\tau} ,\quad n\to \IF.
\EQN
Borrowing the idea of \HL \, we provide next a refinement of \eqref{eq:asT}.\\

\BT \label{theo:2} Let $\X, \X^{(1)} \ldot \X^{(n)}$ be independent Kotz Type III random vectors in $\R^d,d\ge 2$ with distribution function $G$.
Assume that $\X=RA \U$ is  such that the associated random radius $R$ has tail asymptotics given by \eqref{eq:kotz:Fk} with $c,\delta,\tau\in (0,\IF), N\inr$. Assume for simplicity that $\LM=1$ and define its multiplicity $m$ as in \netheo{theo:A}. If $\delta \tau=1$, then we have
\BQN\label{eq:last}
\pb{ \max_{1 \le j  \le n} \norm{\X^{(j)}} \ge b_n^* \quad  i.o.} &= &
 \begin{cases}
 1 & \text{if  }  s \ge (d-m)/2 +N/\tau+1,\\
0 & \text{if   }  s <  (d-m)/2 +N/\tau+1,
\end{cases}
 \EQN
 where $b_n^*:= [\tau (\ln n+ s \ln (\ln n))]^{ 1/\tau}, n>1$.
\ET
\begin{remark} 1.  \
A Kotz Type III random vector $\X$ (see Example 3) is a mean zero Gaussian random vector with covariance matrix  $\SI= A A^\top$ if $N= d-2, \delta=1/2, \tau=2$ implying that \eqref{eq:last} holds with
$$b_n^*:= [ 2 (\ln n + \frac{m}{2} \ln (\ln n))]^{1/2},\quad n>1,$$
 which is shown in Theorem 3 of \HL.

  2. \ Under the Gaussian setup as shown in \HL \, the result of \netheo{theo:2} can be utilised as a diagnostic tool to detect departure from the Gaussian distribution. In view  of our more general result this same tool can be employed for detecting departure from the Kotz Type III multivariate distribution. In the Gaussian case the term $N/\tau$  equals $d/2-1$. In the more general setup of Kotz Type III multivariate distribution $N/\tau$ is in general unknown and needs to be estimated.

 3. \  Further refinements of \eqref{eq:last} dealing also with the case $\delta \tau \not= 1$
can be achieved by borrowing the ideas presented in Embrechts et al.\ (1997); see Example 3.5.6 and Example 3.5.8 therein.
 \end{remark}

\subsection{Asymptotics in the Weibull Model}
Next we consider elliptical random vectors $\X= A R \U$ where $R$ has distribution function  $F$  in the Weibull max-domain of attraction. Necessarily the upper endpoint $x_F$ of $F$ is finite.
Specifically, we suppose that $x_F=1$ and \eqref{eq:Psi} holds for some $\gamma \in  (0, \IF)$.
A canonical example of $F$ in the Weibull max-domain of attraction is the Beta distribution (see Example 4 below).
As we show in the next result, the asymptotic behaviour of $\norm{\X}$ is determined by the eigenvalues of $\SI$ and the tail asymptotics of $F$. \\

\BT \label{theo:B} Under the assumptions and the notation of \netheo{theo:A} if further $x_F= 1$ and $F$ is in the max-domain of attraction of $\Psi_\gamma, \gamma \in (0,\IF)$, then we have $(u \downarrow 0)$
\BQN\label{eq:theo:B:1}
\pk{\norm{\X} > (1- u) \sqrt{ \LM} } 
&=& (1+o(1)) C_*   \frac{\Gamma(\gamma+1)}{\Gamma(\gamma+ (d-m+2)/2)}
\frac{\Gamma(d/2)}{\Gamma(m/2) } (2 u)^{(d-m)/2}\overline{F}(1- u).
\EQN
Furthermore, $\norm{\X}$ possesses the positive density function $h$ with asymptotic behaviour
\BQN\label{eq:theo:A:2}
h(1- u)  &=& (1+o(1)) (\gamma+ (d-m)/2)\pk{\norm{\X} > 1- u}/u, \quad u\downarrow 0.
\EQN
\ET
As in the Gumbel setup we utilise \eqref{eq:theo:B:1} to derive the asymptotics of the sample maxima. 

\BK \label{korr:2} Let $\X, \X^{(1)} \ldot \X^{(n)},n\ge 1$ be independent elliptical random vectors in $\R^d,d\ge 2$ with common distribution function $G$ such that $\X$ satisfies the assumptions of \netheo{theo:B} and  $\LM=1$. Then the convergence in distribution
\BQN\label{eq:conv:dis:B}
\frac{\max_{1 \le j \le n} \norm{\X^{(j)}} - 1}{H^{-1}(1- 1/n)}  &\todis & Y \sim \Psi_{\gamma+(d-m)/2} , \quad n\to \IF
\EQN
holds with  $H^{-1}$ the generalised inverse of the distribution function of $\norm{\X}$.\\
Furthermore, \eqref{eq:conv:dis:B} can be strengthened to the local uniform convergence of the corresponding density functions.
\EK
We give next an illustrating example.

\bigskip
{\bf Example 4.} Let $\X = A R \U$ be an elliptical random vector in $\R^d,d\ge 2$. We consider the special case that
$R \sim Beta(a,b)$ with $a,b$ two positive constants. Since
$$ \pk{R> 1- u}= (1+o(1))\frac{\Gamma(a+b)}{\Gamma(a)\Gamma(b+1)} u^b, \quad u \downarrow 0$$
it follows that $R$ has distribution function in the Weibull max-domain of attraction with index $b$. Consequently,
under the assumptions of \netheo{theo:B} we obtain
\BQNY
\pk{\norm{\X} > (1- u) \sqrt{ \LM} }
&=& (1+o(1)) 2^{(d-m)/2}C_*   \frac{\Gamma(a+b)}{\Gamma(a) \Gamma(b+ (d-m+2)/2)}
\frac{\Gamma(d/2)}{\Gamma(m/2) } u^{(d-m)/2 + b} , \quad u\downarrow  0.
\EQNY

\section{Related Results and Proofs}
\BT \label{eq:theo:BM:A} Let  $R \sim F, X \sim H, Z_{a,b} \sim Beta(a,b),a,b>0$
 be three independent random variables.  Suppose that $x_F\in \{1, \IF\}, F(0)=H(0)=0$ and $H$ has upper endpoint 1.
 Assume further that $H$ is in the max-domain of attraction of $\Psi_{\lambda}, \lambda \in (0,\IF)$ and
 set $Y:= R[(X- \delta) Z_{a,b} + \delta]$ with $\delta\in [0, 1)$.\\
a)  If $F\in MDA(\Lambda, w)$ with some positive scaling function $w$, then we have
\BQN \label{res:BERM1:Aa}
\pk{  Y > u} &=& (1+o(1))
\frac{\Gamma(\lambda+1) \Gamma(a+b)}{ \Gamma(a)} [(1- \delta) u w(u)]^{-b}
\overline{F}(u)\overline{H}(1- 1/(u w(u)) , \quad u \uparrow x_F.
\EQN
Furthermore, the random variable $Y$ possesses  a positive density function $q$ such that
\BQN\label{res:BERM1:Ab}
\lim_{u \uparrow x_F} \frac{q(u)}{ w(u) \pk{Y>u}}&=&1.
\EQN
b) Suppose that  $F$ satisfies \eqref{eq:Psi} with $\gamma\in (0,\IF)$ and $x_F=1$.  Then we have
\begin{eqnarray}
\label{res:BERM1:Ac}
\pk{ Y > u}&=& (1+o(1))
  \frac{\Gamma(a+b)}{\Gamma(a)}
\frac{\Gamma(\lambda+1)\Gamma( \gamma+1)}{\Gamma(\gamma+b+ \lambda+1)}
[(1- \delta)/(1-u)]^{-b} \overline{F}(u)\overline{H}(u), \quad u \uparrow 1.
\end{eqnarray}
Moreover, the density function $q$ satisfies
\BQN\label{eq:BERM1:Ad}
\lim_{u \uparrow 1} \frac{ (1-u) q(u)}{\pk{Y > u}}&=&\gamma+ \lambda+b.
\EQN
\ET
{\bf Proof:} a) We show next the proof only for $x_F=\IF$. When $F$ has a finite upper endpoint the proof follows with similar arguments
utilising further \eqref{eq:uu}, therefore is omitted here. Since $F$ is rapidly varying (see e.g., Resnick (2008)) i.e.,  $\limit{u} \overline{F}(uc)/\overline{F}(u)=0$
for any $c\in (1, \IF)$, and recalling that $R$ is independent of the random variable $(X- \delta) Z_{a,b}+ \delta$, as in the proof of Theorem 12.3.1 in Berman (1992) for any $\ve>0$ we have (set $Y:= R[(X- \delta) Z_{a,b} + \delta]$)
\BQNY
\pk{Y > u} &=& (1+o(1))\int_{u}^{u(1+ \ve)} \pk{(X- \delta) Z_{a,b}+ \delta> u/r} \, d F(r), \quad u\to \IF.
\EQNY
Define  for any $u,s,y\in (0,\IF)$
 $$ \eta(u):= (uw(u))^{-1}, \quad \delta_u:=\eta(u)/(1- \delta),
 \quad H_u(y):= H(1- y \eta(u)), \quad F_u(s):= F(u+ s/w(u)).$$
Transforming the variables we have
\BQNY
\pk{Y > u} &=& \int_{0}^{\ve / \eta(u)} \int_{0}^{s+ o(1)}  \pk{Z_{a,b}> ((1+ s \eta(u))^{-1}  - \delta)/(1- y \eta(u) - \delta) }  \, d H_u(y)\, d F_u(s)\\
 &=& \int_{0}^{\ve / \eta(u)} \int_{0}^{s+o(1)} \pb{ Z_{a,b} > 1- \delta_u(s- y)(1+o(1)) }  \, d H_u(y) \, d F_u(s).
 \EQNY
By the assumptions on $F$ and $H$ we may write
$$ \lim_{u \uparrow x_F}  \frac{F_u(s) - F_u(t)}{\overline{F}(u)}= \exp(-t)- \exp(-s), \quad  \lim_{u \uparrow x_F}  \frac{ H_u(s) - H_u(t)}{\overline{H}(u)}= t^\lambda- s^\lambda, \quad  \forall s,t\inr.$$
Furthermore, since $Z_{a,b}\sim Beta(a,b)$ for any $s,y\inr$ such that $s\ge y$ we have (recall  \eqref{eq:uu})
$$ \lim_{u \uparrow x_F} \delta_u^{-b} \pb{ Z_{a,b} > 1- \delta_u(s- y)(1+o(1)) }  = c_{a,b}(s-y)^b, $$
with $ c_{a,b}:= \Gamma(a+b)/(\Gamma(a)\Gamma(b+1)).$ Hence Fatou Lemma implies
\BQNY
\liminf_{u \uparrow x_F}
\frac{ \pk{Y > u}}{ \delta_u^{b} \overline{F}(u)\overline{H}(1- \eta(u))}
&\ge &
c_{a,b}\lambda \int_{0}^{\IF} \exp(-s) \int_{0}^s (s-y)^b  y^{\lambda -1} \, dy \, ds \\
&=&   c_{a,b}\lambda \int_{0}^{\IF} \exp(-s) s^{b+ \lambda} \int_{0}^1 (1-y)^b  y^{\lambda -1} \, dy \, ds \\
&=& c_{a,b}\lambda  \Gamma(\lambda+b+1)\frac{ \Gamma(b+1) \Gamma(\lambda)}
{\Gamma(\lambda+b+1)}\\
&=& \frac{\Gamma(\lambda+1)\Gamma(a+b)}{\Gamma(a)}.
 \EQNY
The proof for the $\limsup$ (which coincides with $\liminf$) can be established
utilising Lemma 7.5 an Lemma 7.7 in Hashorva (2007a).\\
Next, the independence of $X$ and $Z_{a,b}$ and the fact that $Z_{a,b}$ possesses a positive density function in $(0,1)$ implies that
the random variable $X^*:= (X- \delta)Z_{a,b}+ \delta$ possesses a positive density function in $(0,1)$ which we denote by $g$.
Consequently, since $X^*$ and $R$ are independent and $X^*$ possesses the density function $g$ it follows that $Y$ possesses the density function $q$ given by
\BQN\label{eq:den:yy}
 q(u)&=& \int_{u}^{x_F} g( u/r) \frac{1}{r} \, d F(r), \quad u\in (0,x_F).
\EQN
The asymptotic behaviour of $q(u), u\uparrow x_F$ can be establishes with similar arguments as the proof above leading thus
to \eqref{res:BERM1:Ab}.\\
b) Define next for any $u\in (0,1)$
 $$ H_u(y):= H(1- yu), \quad F_u(s):= F(1- su), \quad s,y\in (0,\IF), \quad \delta_u:=u/(1-\delta).$$
We may further write ($u \downarrow 0)$
\BQNY
\pk{Y > 1-u}
 &=& \int_{0}^{1} \int_{0}^{1-s+o(1)} \pb{ Z_{a,b} > 1- \delta_u(1-s- y)(1+o(1)) }  \, d H_u(y) \, d F_u(s).
 \EQNY
By the assumptions on $F$ and $H$ we have
$$ \lim_{u \downarrow 0}  \frac{ F_u(s) - F_u(t)}{\overline{F}(u)}= t^{\gamma}- s^{\gamma}, \quad \lim_{u \downarrow 0}  \frac{ H_u(s) - H_u(t)}{\overline{H}(u)}= t^\lambda- s^\lambda, \quad \forall s,t\inr.$$
As above for any $s,y\inr$ such that $s+y \le 1$
$$ \lim_{u \downarrow 0} \delta_u ^{-b} \pb{ Z_{a,b} > 1- \delta_u (1-s- y)(1+o(1)) }  =c_{a,b}(1-s-y)^b .$$
Hence we obtain applying Lemma 4.2 in Hashorva (2007b)
\BQNY
\frac{\pk{Y > 1-u}}{\delta_u^{b}\overline{F}(1-u)\overline{H}(1-u) } &=& (1+o(1)) c_{a,b}\lambda \gamma
 \int_{0}^{1} s^{\gamma -1}\int_{0}^{1-s} (1-s- y)^b y^{\lambda- 1}  \, d y\, d s \\
 &=& (1+o(1)) c_{a,b} \lambda \gamma
 \frac{\Gamma(b+1)\Gamma(\lambda)}{\Gamma(b+ \lambda+1)} \int_{0}^{1} s^{\gamma -1}(1-s)^{b + \lambda} \, ds \\
 &=& (1+o(1))c_{a,b} \lambda \gamma
\frac{\Gamma(b+1)\Gamma(\lambda)}{\Gamma(b+ \lambda+1)} \frac{ \Gamma(b+ \lambda+1)\Gamma(\gamma)}{\Gamma(\gamma+b+ \lambda+1)}\\
 &=& (1+o(1))   \frac{\Gamma(a+b)}{\Gamma(a)}
\frac{\Gamma(\lambda+1)\Gamma( \gamma+1)}{\Gamma(\gamma+b+ \lambda+1)}, \quad u \downarrow 0.
\EQNY
The proof of \eqref{eq:BERM1:Ad} follows with similar arguments utilising further \eqref{eq:den:yy}.\QED

In the next theorem we derive the asymptotic tail behaviour of the product $X Z_{a,b}$.
Its proof is similar to
that of Theorem 12.3.1 of Berman (1992) (also Berman (1982, 1983), and Hashorva (2007a)), therefore we omit it here.
See Tang  and Tsitsiashvili  (2004) and Tang (2006,2008) for recent results on the tail asymptotics of products of random variables.\\

\BT \label{eq:theo:BM1} Let  $X \sim F, Z_{a,b}\sim Beta(a,b) ,a,b>0$
 be two independent random variables and let $Y := X [1- \delta Z_{a,b}]^{\tau}, \delta\in (0,1], \tau\in (0,\IF)$.
  Suppose that $x_F\in \{1, \IF\}$ and $F(0)=0$.\\
a)  If $F\in MDA(\Lambda, w)$ with some positive scaling function $w$, then  we have
\BQN \label{res:BERM1:a}
\pk{  Y> u} &=& (1+o(1)) \frac{\Gamma(a+b)}{\Gamma(b)} (  \delta \tau uw(u))^{-a} \overline{F}(u), \quad u \uparrow x_F.
\EQN
Furthermore, the random variable $Y$ possesses  a positive density function $q$ and
\BQN\label{res:BERM1:b}
\lim_{u \uparrow x_F} \frac{q(u)}{ w(u) \pk{Y>u}}&=&1.
\EQN
b) Suppose that  $F$ satisfies \eqref{eq:Psi} with $\gamma\in (0,\IF)$ and $x_F=1$.  Then we have
\begin{eqnarray}
\label{res:BERM1:c}
\pk{ Y > u}&=& (1+o(1))\frac{\Gamma(\gamma+1)\Gamma(a+b)}{\Gamma(b)\Gamma(\gamma+a+1)}
((1-u)/(\tau \delta) )^a \overline{F}(u), \quad u\uparrow 1.
\end{eqnarray}
Moreover, the density function $q$ satisfies
\BQN\label{eq:T2:2}
\lim_{u\uparrow 1} \frac{ (1-u) q(u)}{\pk{Y > u}}&=&\gamma+ a.
\EQN
\ET
\COM{
{\bf Proof:} a) The proof of \eqref{res:BERM1:a} can be established along the lines of the proof of Theorem 8 in Hashorva (2007a). See also the proof of Theorem 12.3.1 in Berman (1992). \\
Next, in order to show \eqref{res:BERM1:b} note first that $Y$ possesses a density function  $h$ given by
\BQN\label{eq:dens:h}
 h(s)&=& \int_s^{x_F} g(s/t) \frac{1}{t}\, d F(t), \quad \forall s\in (0,x_F),
 \EQN
with $g$ the positive density function of $Y \equaldis [1- \delta Z_{a,b}]^{\tau}$.
Suppose that $x_F=\IF$. By the assumptions on $F$ we have $\limit{u} u w(u)= \IF$.
Transforming the variables we have (set $v(u):=u w(u), F_u(t):= F(u+  t/w(u)),t\inr, u>0, \overline{F}:= 1- F$)
\BQNY
h(u)&= &\int_u^{\IF } g(u/t) \frac{1}{t}\, d F(t)\\
&= &\int_0^{\IF} g(u/(u+ t/w(u))) \frac{1}{u+ s/w(u)}\, d F_u(t)\\
&= & \frac{1}{u}\overline{F}(u)  g(1/(1+ 1/v(u)))
\int_0^{\IF} \frac{g(1/(1+ t/v(u)))}{
 g(1/(1+ 1/v(u)))} \frac{1}{1+ s/v(u)}\, d F_u(t)/\overline{F}(u), \quad u>0.
\EQNY
Since  as  $u\to \IF$
$$ \frac{g(1/(1+ t/v(u)))}{ g(1/(1+ 1/v(u)))} \to t^{a-1}, \quad \forall  t>0$$
and $F\in MDA(\Lambda, w)$ implies
$$ \frac{F_u(t)- F_u(s)}{\overline{F}(u)}= \exp(-s) - \exp(-t), \quad \forall s,t\inr,$$
then  applying Fatou Lemma we obtain
 $$\liminf_{u\to \IF} \int_0^{\IF} \frac{g(1/(1+ t/v(u)))}{
 g(1/(1+ 1/v(u))) }\, d F_u(t)/\overline{F}(u)\ge
\int_0^{\IF} x^{a-1} \exp(- x) \, dx= \Gamma(a)\in (0,\IF).$$
Consequently
\BQNY
h(u)&= & (1+o(1)) \frac{\Gamma(a+b)}{\Gamma(b)} (  \delta \tau u)^{-a} (w(u))^{1- a}\overline{F}(u), \quad u \uparrow x_F
\EQNY
holds.  By the properties of the density function $g$
utilising Lemma 7.5 an Lemma 7.7 in Hashorva (2007a) it follows that
 $$\limsup_{u\to \IF} \int_0^{\IF} \frac{g(1/(1+ t/v(u)))}{
 g(1/(1+ 1/v(u))) }\frac{1}{1+ t/v(u)}\, d F_u(t)/\overline{F}(u)\le \Gamma(a)\in (0,\IF).$$
The proof for $x_F\in (0,\IF)$ follows with similar arguments making further use of the fact that $\lim_{u\uparrow x_F} w(u) (x_F- u) = \IF$. Thus the first claim follows.\\
b) For $\delta=1$ the result in \eqref{res:BERM1:b} can be found in Theorem 8 in Hashorva (2007a). The proof for $\delta\in (0,1)$
follows with similar argument as that of the aforementioned theorem. We show next the last claim. In view of \eqref{eq:dens:h}
for any $u\in (0,1)$ we have (define $F_u(x):=F(1- ux), u\in (0,1), x\in (0,\IF)$)
\BQNY
\frac{ h(1-u)}{\overline{F}(1-u)}
&=& \int_{1- u}^{1} g((1-u)/t) \frac{1}{t}\, d F(t)/\overline{F}(1-u)\\
&=& \int_{0}^{1} g(1-u(1- x)(1+o(1))) \frac{1}{1- ux}\, d F_u(x)/\overline{F}(1-u).\\
 \EQNY
As $u \downarrow 0$
$$ \frac{g(1- u(1- x)(1+o(1)))}{ g(1- u} \to (1-x)^{a-1}, \quad \forall  x\in (0,1),$$
thus applying Lemma 4.2 in Hashorva (2007b) we obtain
\BQNY
\frac{ h(1-u) g(1- u)}{\overline{F}(1-u)}
&=& (1+o(1))\int_{0}^{1} (1- x)^{a+\gamma-1} \, dx , \quad u \downarrow 0,
 \EQNY
hence the result follows easily . \QED
}

\prooftheo{theo:A} If the multiplicity of $\LM$ is $d$, i.e., $\lambda_1 = \lambda_2 =\cdots =\lambda_d$, then 
$$\norm{\X}^2=  \LM R^2 \norm{\U}^2= \LM R^2,$$
hence the claim follows. Assume next that $m< d$, and define a random vector $\vk{V}_K$ such that
$\vk{V}_K\equaldis  R \sqrt{W_{m, d}} \Um,$
with $R\sim F$ being independent of the random variable $W_{m, d} \sim Beta(m/2, (d-m)/2)$.  We have
$$ \norm{\vk{V}_K}^2 \equaldis   R^2 W_{m,d} \norm{\Um}^2\equaldis   R^2 W_{m,d}= R_m^2,$$
with $R_m:= R \sqrt{W_{m, d}}.$ Consequently since $R$ and $W_{m,d}$ are independent applying \netheo{eq:theo:BM1} we obtain
$$\pk{ \norm{\vk{V}_K} > u}=
(1 + o(1)) \frac{\Gamma(d/2)}{  \Gamma(m/2)}  \fracl{2}{ u w(u) }^{ (m-d)/2} \overline{F}(u), \quad u\uparrow x_F.$$
By the self-neglecting property (recall \eqref{eq:self}) of the scaling function $w$
we conclude that the random variable $\norm{\vk{V}_K}$ has distribution function in
the Gumbel max-domain of attraction with the same scaling function $w$.

Next, since $\X \equaldis AR\U $ is an elliptical random vector we may write
$$ \norm{\X}^2 \equaldis \LM \norm{\Y}^2,$$
where $\Y \equaldis  D R \U$ with $D$ a diagonal matrix with $d_{ii}= 1, i=1 \ldot m$ and $d_{ii}=\sqrt{\lambda_i/\LM}, i=m+1 \ldot d$. Note that $d_{ii} \in (0,1)$ for all $i> m$. Define $K_i:= K_{i-1} \cup \{i\}, i:= m+1 \ldot d,$ with $K_m:=K$.
In the light of \eqref{rep:elli0}
$$\Y_{K_{m+1}} \equaldis R_m \Bigl( \Um \sqrt{W_{m, m+1}}, I_1 \sqrt{(1- W_{m,m+1})\lambda_{m+1}/\LM}\Bigr),$$
 where $I_1$ assumes only two values $-1,1$ with equal probability, $W_{m,m+1}\sim Beta(m/2,1/2)$ and the random variables
 $I_1, R_m, \Um,W_{m,m+1}$ are mutually independent.  Consequently,
$$ \norm{\vk{X}_{K_{m+1}}}^2=  \LM R^2_m[ W_{m, m+1}+  \lambda_{m+1} /\LM (1- W_{m,m+1})]\equaldis
\LM R^2_m[ 1- (1- \lambda_{m+1} /\LM) W_{ m,m+1}^*]$$
holds with $W_{m,m+1}^* \sim Beta(1/2, m/2)$ being independent of $R_m$. Hence applying again \netheo{eq:theo:BM1} we obtain
$$ \pk{ \norm{\X_{K_{m+1}}} >u \sqrt{\LM } }=(1+o(1)) (1- \lambda_{m+1}/\LM)^{-1/2} \pk{\norm{\vk{V}_K}> u}, \quad u\uparrow x_F.$$
Similarly 
$$\norm{\X_{K_{m+2}}}^2 \equaldis (\LM R^2_m[ 1- (1- \lambda_{m+1} /\LM) W_{ m,m+1}^*])W_{ m+1,m+2}^* + \lambda_{m+2} (1- W_{ m+1,m+2}^*),$$
 where $R_m, W_{ m,m+1}^*, W_{ m+1,m+2}^*$ are independent and $W_{ m+1,m+2}^*\sim Beta((m+1)/2,1/2)$. Applying \netheo{eq:theo:BM:A}
we have
$$ \pk{  \norm{\X_{K_{m+2}}} >u \sqrt{\LM }}=(1+o(1)) (1- \lambda_{m+1}/\LM)^{-1/2}
(1- \lambda_{m+2}/\LM)^{-1/2} \pk{\norm{\vk{V}_K}> u}, \quad u\uparrow x_F.$$
The proof follows now by applying \netheo{eq:theo:BM1} iteratively and utilising \eqref{rep:elli0}.
\QED

\proofkorr{korr:1} The proof follows from the result of \netheo{theo:A} and the self-neglecting property of $w$ in \eqref{eq:self}.
\QED

\prooftheo{theo:2} The proof follows by the formula \eqref{tail:K:III} utilising further Lemma 1 in \HL. \QED

\prooftheo{theo:B} The proof follows along the lines of the proof of \netheo{theo:A} utilising further \netheo{eq:theo:BM1} and making use of
the asymptotic condition \eqref{eq:de:W}. \QED

{\bf Acknowledgment:} I would like to thank Andre Moser and Professor J\"urg H\"usler for several motivating discussions. To both referees of the paper I am in debt for kind review containing several suggestions and corrections.

\bibliographystyle{plain}

\end{document}